\documentclass[english]{article}
\usepackage{amsmath,amsthm,eucal,upref,cite,calc,graphicx,tabularx,hhline,longtable,enumerate,comment}
\newtheorem{thrm}{Theorem}
\newtheorem{lem}{Lemma}
\usepackage{amssymb}
\usepackage{hyperref}
\theoremstyle{plain}
\newcommand{\al}{\alpha}
\newcommand{\be}{\beta}

\newcommand{\ep}{\varepsilon}

\newcommand{\M}{\mathcal M}

\renewcommand{\le}{\leqslant}
\renewcommand{\ge}{\geqslant}

\newcommand{\const}{\mathrm{const}}

\renewcommand{\atop}[2]{\genfrac{}{}{0pt}{1}{#1}{#2}}

\renewcommand{\P}{\mathsf{P}}
\newcommand{\E}{\mathsf{E}}

\newcommand{\pf}{\begin{proof}}
\newcommand{\epf}{\end{proof}}

\begin{document}

\title{On some estimates for Erd\"os-R\'enyi random graph}
\author{Nikolay Kazimirow}
\date{8/31/2015}
\maketitle

\begin{abstract}
We consider a number $\nu_n$ of components in a random graph $G(n,p)$ with $n$ vertices, where the probability of an edge is equal to $p$.
By operating with special generating functions we shows the next asymptotic relation for factorial moments of $\nu_n$:
$$
\E(\nu_n-1)^{\underline s} = (1+o(1))\left( \frac 1p \sum\limits_{k=1}^\infty\frac{k^{k-2}}{k!}(npq^n)^k\right)^s + o(1)
$$
as $n$ tends to $\infty$ and $q=1-p$. And the following inequations hold:
$$
1-2nq^{n-1} \le p_n\le\frac{1}{nq^n},
$$
$$
1-\frac{1}{nq^n}\le pi_n\le nq^{n-1},
$$
where $p_n$ is the probability that $G(n,p)$ is connected and $pi_n$ is the probability that $G(n,p)$ has an isolated vertex.
\end{abstract}


\section{Notations}

Let $G_n$ be a set of undirected graphs with $n$ labeled vertices.
For any graph $g\in G_n$ let $C(g)$ be a number of connected components in the graph $g$ and
$E(g)$ be a number of edges in thr graph $g$. Besides we denote by $F_{s,n}$ the number of all forests in $G_n$,
that contains exactly $s$ trees. We also suppose that components in $G_n$ are not ordered.

Further, let $A_{n,k,s}$ be a number of graphs in $G_n$, which contains $n$ vertices, $k$
edges and $s$ components, $A_{n,k}$ be a number of graphs, which contains $n$
vertices and $k$ edges, and $B_{n,k}$ --- a number of connected graphs with
$n$ vertices and $k$ edges. For definiteness we suppose that
$A_{0,k}=A_{0,k,s}=A_{n,k,0}=0$ in all cases, except $n=k=s=0$, where we set by definition
$A_{0,0}=A_{0,0,0}=1$. Besides, let $B_{0,k}=0$ for all $k$. It's clear that $A_{n,k}=\sum_s A_{n,k,s}$,
where index $s$ runs on all integer non-nagative numbers.

Let us consider the random graph $G(n,p)$, which contains $n$ labeled vertices, where each of $\binom n2$ edges is present with the probability $p$ independently of other edges. Each concrete realization of random graph $G(n,p)$ is a graph from $G_n$. This model of random graphs was firstly described by Erd\"os and R\'enyi in \cite{Erdos-Renyi,Erdos-Renyi2} and then has been well studied by B\'ela Bollob\'as \cite{Bollobas}, Valentin Kolchin \cite{Kolchin} and other authors.

It is easy to see that the parobability distribution of such random graph is defined as follows:
$$
\P\{G(n,p)=g\}=(p/q)^{E(g)}q^{n(n-1)/2},
$$
where $g\in G_n$ and $q=1-p$.

Let denote by $\nu_n$ the number of connected components of $G(n,p)$, i.\,e. $\nu_n=C(G(n,p))$, and let $p_n$ be the probability that random graph $G(n,p)$ is connected, thus
$p_n=\P\{\nu_n=1\}$. It's clear that
\begin{gather}
\P\{\nu_n=s\}=\sum_{k=0}^\infty A_{n,k,s}(p/q)^{k}q^{n(n-1)/2}\label{nus}
\end{gather}
and
$$
p_n=\sum_{k=0}^\infty B_{n,k}(p/q)^{k}q^{n(n-1)/2}.
$$
Froom the above agreements it follows that $p_0=0$ and $p_1=1$.

Below we'll need the special generated function, which we define as follows: for a sequence of functions $\{r_n(q)\}$ we put
$$
R=R(x,q)=\sum_{n=0}^\infty \frac{x^n}{n!q^{n(n-1)/2}}r_n(q),
$$
where we often will skip arguments $x$ and $q$, except such cases when we will use special values of them. Below in this text we will call such functions
as SG-functions (SG = special generated).

It is easy to see that SG-functions are formal power series which are not converges at all. But most of all usual operations with SG-functions (such as adding, production, differentiation and integration on both arguments) does not lead to conflicts when counting coefficients before $x^n$.

Let denote
$$
\widehat R =\sum_{n=0}^\infty\frac{x^n}{n!q^{n(n-1)/2}}\frac{dr_n(q)}{dq},
$$
i.e. the operator $\widehat{}$ denotes SG-function for the sequence of derivatives of $r_n(q)$.

Let also:
\begin{align*}
A&=\sum_{n=0}^\infty\frac{x^n}{n!q^{n(n-1)/2}}\\
B&=\sum_{n=0}^\infty\frac{x^n}{n!q^{n(n-1)/2}}p_n\\
E&=\sum_{n=0}^\infty\frac{x^n}{n!q^{n(n-1)/2}}\E\nu_n\\
M_k&=\sum_{n=0}^\infty\frac{x^n}{n!q^{n(n-1)/2}}\E(\nu_n)^{\underline k}\\
\M_k&=\sum_{n=0}^\infty\frac{x^n}{n!q^{n(n-1)/2}}\E(\nu_n-1)^{\underline k},
\end{align*}
where $z^{\underline k}=z(z-1)\dots(z-k+1)$ denotes the factorial power $k\ge 0$. Therefore,
$A$ is a SG-function of $\{1\}$, $B$ is a SG-function of probabilies that graph is connected, $E$ is a SG-function of expectations of components quantity, $M_k$ is a SG-function of $k$-th factorial moments of $\nu_n$, and $\M_k$ is a SG-function of $k$-th factorial moments of $(\nu_n-1)$. It's easy to see that $A$ converges only if $q=1$. Below we'll see that all of theese series are converges in the same conditions.

\section{Basical Relations}

\begin{lem}
If the relation $n=k=s=0$ does not holds, then
\begin{equation}\label{A-B}
A_{n,k,s}=\sum_{\atop{n_1+\dots+n_s=n}{k_1+\dots+k_s=k}}\frac{n!}{s!}\frac{B_{n_1,k_1}\cdots B_{n_s,k_s}}{n_1!\cdots n_s!},
\end{equation}
where the summation is over all integer non-negative $n_i$, $k_i$.
\end{lem}
\pf
Let consider the set of graphs $\bar G_n$ with $n$ vertices, where the components are ordered.
It is clear that the number $\bar A_{n,k,s}$ of such graphs with $n$ vertices, $k$ edges and $s$ components is equal to $s!A_{n,k,s}$.

By the other side, any graph from $\bar G_n$ with $n$ vertices and $s$ components we can make by getting some ordered partition of the set of $n$ vertices with non-empty parts, which has the volumes $n_1,\dots,n_s$. The number of such partitions is equal to $n!/(n_1!\cdots n_s!)$. For every set of vertices, included in connected components, we can find the number of connected graphs with $n_i$ vertices and $k_i$ edges. It is equal to $B_{n_i,k_i}$.
By choosing $k_i$ in such a way that $k_1+\dots+k_s=k$, and summing over all partitions of $n$ vertices, we get the equation:
$$
\bar A_{n,k,s}=\sum_{\atop{n_1+\dots+n_s=n}{k_1+\dots+k_s=k}}\frac{n!B_{n_1,k_1}\cdots B_{n_s,k_s}}{n_1!\cdots n_s!}
$$
From this we get \eqref{A-B} for positive $n,n_i,s$ and non-negative $k$. Extention of this relation for zero values of $n,n_i$ and $s$
follows from the previous agreements.
\epf

Now we consider the next generated functions, which are exponential by parameter $x$:
$$
A(x,y)=\sum_{n,k}\frac{A_{n,k}}{n!}x^ny^k,\qquad B(x,y)=\sum_{n,k}\frac{B_{n,k}}{n!}x^ny^k.
$$
The summation is over integer non-negative $n,k$.
\begin{lem}
\begin{equation}\label{pf}
A(x,y)=e^{B(x,y)}
\end{equation}
\end{lem}
\pf By multyplying the relation \eqref{A-B} by $x^ny^k/n!$ we get:
$$
\frac{A_{n,k,s}}{n!}x^ny^k=%
\frac{1}{s!}\sum_{\atop{n_1+\dots+n_s=n}{k_1+\dots+k_s=k}}\frac{B_{n_1,k_1}x^{n_1}y^{k_1}\cdots B_{n_s,k_s}x^{n_s}y^{k_s}}{n_1!\cdots n_s!}=%
[x^ny^k]B(x,y)^s.
$$
The last notation denotes a coefficient before $x^ny^k$ in the series $B(x,y)^s$. Now, by summing over integer non-negative
$n,k$ for $s>0$ we get the following:
\begin{equation}\label{B^s}
\sum_{n,k}\frac{A_{n,k,s}}{n!}x^ny^k=\frac{1}{s!}B(x,y)^s
\end{equation}

Note, that by virtue of the agreements this equation stays also true for $s=0$.
Finally, by summing over integer non-negative $s$ we get:
$$
A(x,y)=\sum_{s=0}^\infty\frac{B(x,y)^s}{s!}=e^{B(x,y)}.
$$
\epf

From the relation \eqref{pf} we can obtain any exact expressions for probabilities of random graph $G(n,p)$. First of all, it is clear that:
$$
A_{n,k}=\binom{n(n-1)/2}{k},
$$
where we suppose that $\binom mk=0$ for $k>m$. It is easy to see that
$$
\sum_{k=0}^\infty \binom{n(n-1)/2}{k}y^k=\sum_{k=0}^{n(n-1)/2}\binom{n(n-1)/2}{k}y^k=(1+y)^{n(n-1)/2},
$$
hence,
$$
A(x,y)=\sum_{n=0}^\infty (1+y)^{n(n-1)/2}\frac{x^n}{n!}.
$$

From this and from \eqref{pf} it follows that
\begin{equation}\label{pf-sv}
B(x,y)=\ln\sum_{n=0}^\infty (1+y)^{n(n-1)/2}\frac{x^n}{n!}.
\end{equation}

One can see that $B(x,y)$ is the generated function for a sequence \cite{A062734} where the nulled element is equal to zero.

By putting $y=p/q$ and from the obvious equtions $$\sum\limits_k A_{n,k}\left(\frac pq\right)^kq^{n(n-1)/2}=1,\quad 
\sum\limits_k B_{n,k}\left(\frac pq\right)^kq^{n(n-1)/2}=p_n,$$ we get that for previously defined series $A$ and $B$
the next relations are true:
\begin{gather}
A\left(x,\frac pq\right)=\sum_{n,k}\frac{A_{n,k}}{n!}x^n(p/q)^k=\sum_n\frac{x^n}{q^{n(n-1)/2}n!}=A,\nonumber\\[-6pt]
\label{Bx}\\[-6pt]
B\left(x,\frac pq\right)=\sum_{n,k}\frac{B_{n,k}}{n!}x^n(p/q)^k=\sum_n\frac{p_nx^n}{q^{n(n-1)/2}n!}=B.\nonumber
\end{gather}

Thus, we have
\begin{lem}\label{AeB}
$$
A=e^B.
$$
\end{lem}
\noindent This proved equation is the base fact, which we will use anythere below without a special link.

From \eqref{nus} it follows that:
$$
\sum_{n=0}^\infty\frac{\P\{\nu_n=s\}}{q^{n(n-1)/2}}\frac{x^n}{n!}=\sum_{n,k}\frac{A_{n,k,s}}{n!}x^n(p/q)^k,
$$
and by \eqref{B^s}, where we put $y=p/q$, we get following:
\begin{equation}\label{Bs}
\sum_{n=0}^\infty\frac{x^n}{q^{n(n-1)/2}n!}\P\{\nu_n=s\}=\frac{1}{s!}B(x,p/q)^s=\frac{1}{s!}B^s,
\end{equation}
i.e. the formal series $B^s/s!$ is SG-function of probabilities $\P\{\nu_n=s\}$ for a fixed number $s$ of connected components.

Let us consider two SG-functions and their product:
$$
R=\sum_{n=0}^\infty\frac{r_nx^n}{n!q^{n(n-1)/2}},\quad 
T=\sum_{n=0}^\infty\frac{t_nx^n}{n!q^{n(n-1)/2}},\quad
RT=\sum_{n=0}^\infty\frac{z_nx^n}{n!q^{n(n-1)/2}}.
$$

One can easily proof the following
\begin{lem}[Convolution Formula]\label{convolution} For $n\ge 0$:
$$
z_n=\sum_{k=0}^n\binom nk q^{k(n-k)}r_kt_{n-k}.
$$
\end{lem}
\noindent Further we will use this formula without a special link to it. The next recursion formula for probabilities $p_n$ is an anlogue of a recursion formula for a number of connected graphs, that was obtained in \cite{Harari-Palmer}.
\begin{lem} For any $n\ge 1$
\begin{equation}\label{pn-rekurs}
p_n=1-\sum_{k=1}^{n-1}\binom {n-1}{k} q^{k(n-k)}p_{n-k}.
\end{equation}
\end{lem}
\pf
By differentiating the relation $A=e^B$ by the parameter $x$ we get:
$$
xA'=xAB',
$$
hence, from the \hyperref[convolution]{convolution formula} it follows that
\begin{equation}\label{nkpk}
n=\sum_{k=0}^n\binom nk q^{k(n-k)}kp_k
\end{equation}
Since $p_0=0$ and $\frac kn\binom nk=\binom{n-1}{k-1}$ follows
$$
1-p_n=\sum_{k=1}^{n-1}\binom {n-1}{k-1}q^{k(n-k)}p_k,
$$
and by replacing $k$ by $n-k$ we get the statement of Lemma.
\epf

By analogue we can get a recursive formula for probabilities $\P\{\nu_n=s\}$.
\begin{lem}
\begin{equation}\label{rekurs}
\P\{\nu_n=s\}=\sum_{k=s-1}^{n-1}\binom{n-1}k\P\{\nu_k=s-1\}p_{n-k}q^{k(n-k)}
\end{equation}
for $n\ge s>1$.
\end{lem}
\pf Let us denote
$$
B_s=\sum_{n=0}^\infty\frac{x^n}{q^{n(n-1)/2}n!}\P\{\nu_n=s\},
$$
then by \eqref{Bs} we get:
$$
s!B_s(x)=B(x)^s,
$$
then by differentiating by $x$ it follows that:
$$
s!B'_s=sB^{s-1}B'=s(s-1)!B_{s-1}B',
$$
hence,
$$
xB'_s=xB_{s-1}B'.
$$
From this and according to $\P\{\nu_k=s-1\}=0$ as $k<s-1$ we get Lemma statement.
\epf

\begin{lem} The following relations hold:
\begin{gather}\label{pn+1p}
p_{n+1}=\sum_{s=1}^n\sum_{k_1+\dots+k_s=n}\frac{n!(1-q^{k_1})\dots(1-q^{k_s})}{s!k_1!\dots k_s!}\P\{\nu_n=s\},
\end{gather}
\begin{gather}\label{pn+1pn}
p_{n+1}\ge(1-q^n)p_n.
\end{gather}
\end{lem}
\pf If we put $x/q$ instead of $x$ in the definition of series $A$, we get that $A'=A(x/q)=e^{B(x/q)}$.
On the other side, $A'=B'e^B$. Therefore,
$$
B'e^B=e^{B(x/q)},\qquad B'=e^{B(x/q)-B(x)},
$$
hence,
$$
B'=\sum_{s=0}^\infty\frac 1{s!}\left(\sum_{n=0}^\infty\frac{p_nx^n(1-q^n)}{n!q^{n(n-1)/2}}\right)^s.
$$
Now we take the corresponding coefficients before $x^n$ in theese series and get the relation \eqref{pn+1p}.
The ineqution \eqref{pn+1pn} follows from \eqref{pn+1p} if we left in this summa only the summand with $s=1$.
\epf


\section{Several Equations}

\begin{lem} For $s\ge 0$
$$
M_s=AB^s,
$$
and in particulary, $E=AB$.
\end{lem}
\pf
By definition,
$$
\E(\nu_n)^{\underline s}=\sum_{k=0}^nk^{\underline s}\P\{\nu_n=k\},
$$
hence by \eqref{Bs} we get:
\begin{align*}
M_s&=\sum_{n=0}^\infty\frac{\E(\nu_n)^{\underline s}x^n}{n!q^{n(n-1)/2}}=
     \sum_{k=0}^\infty k^{\underline s}\sum_{n=0}^\infty\frac{\P\{\nu_n=k\}x^n}{n!q^{n(n-1)/2}}=\\
   &=\sum_{k=0}^\infty k^{\underline s} B^s/s!=B^s\sum_{k=s}^\infty\frac{B^{k-s}}{(k-s)!}=AB^s.
\end{align*}
\epf

Now we consider the connection between moments of $\nu_n$ and $\nu_n-1$.
\begin{lem}\label{MkMkMk} For $s\ge 1$
$$
M_s=\M_s+s\M_{s-1}
$$
$$
\frac{(-1)^s}{s!}\M_s=\sum_{k=0}^s\frac{(-1)^k}{k!}M_k
$$
\end{lem}
\pf
The first equation is follows from
$$
\E(\nu_n-1)^{\underline s}=\E(\nu_n-1)\dots(\nu_n-s)=\E(\nu_n)^{\underline s}-s\E(\nu_n-1)^{\underline {s-1}},
$$
and the second one not hard to proof by induction with the obvious start eqation $\M_0=A=M_0$.
\epf
\begin{lem}\label{M'kMkMk} For $s\ge 1$
\begin{align*}
\frac{M_s'}{s!}&=B'\left(\frac{M_s}{s!}+\frac{M_{s-1}}{(s-1)!}\right)\\
\frac{\M_s'}{s!}&=B'\left(\frac{\M_s}{s!}+\frac{\M_{s-1}}{(s-1)!}\right)=B'\frac{M_s}{s!}
\end{align*}
\end{lem}

Now we ready to use the operator $\widehat{}$ for SG-functions of moments. First of all, we get:
\begin{lem}[Derivative Relashionship Formula] If $R$ is a SG-function, then:
$$
\widehat R=R'_q+\frac{x^2}{2q}R''
$$
\end{lem}
Here and below the single quote without a parameter notation denotes the derivative by $x$, and the derivative by $q$ is marked by index $q$.

The following equations hold.
\begin{lem}\label{hatM-Mu}
\begin{align*}
\frac{\widehat M_s}{s!}&=\frac{x^2}{2q}(B')^2\left(\frac{M_{s-1}}{(s-1)!}+\frac{ M_{s-2}}{(s-2)!}\right)=\frac{x^2}{2q}B'\frac{M'_{s-1}}{(s-1)!}\\
\frac{\widehat \M_s}{s!}&=\frac{x^2}{2q}(B')^2\left(\frac{ \M_{s-1}}{(s-1)!}+\frac{ \M_{s-2}}{(s-2)!}\right)=\frac{x^2}{2q}(B')^2\frac{M_{s-1}}{(s-1)!}=\frac{x^2}{2q}B'\frac{\M'_{s-1}}{(s-1)!}
\end{align*}
\end{lem}
\pf
By the \hyperref[convolution]{convolution formula} and from $\widehat A=0$ we get:
$$
A'_q=-\frac{x^2}{2q}A''.
$$
From here it follows that:
\begin{align*}
(M_s)'_q&=(AB^s)'_q=A'_qB^s+sAB^{s-1}B'_q=A'_q(B^s+sB^{s-1})=-\frac{x^2}{2q}A''(B^s+sB^{s-1})\\
M_s''&=(AB^s)''=(A'B^s+sB^{s-1}A')'=A''(B^s+sB^{s-1})+A'B'(sB^{s-1}+s(s-1)B^{s-2})
\end{align*}
Hence by Derivative Relashionship Formula we get that:
\begin{align*}
\widehat M_s&=(M_s)'_q+\frac{x^2}{2q}M_s''=\frac{x^2}{2q}A'B'(sB^{s-1}+s(s-1)B^{s-2})=\\
&=\frac{x^2}{2q}(B')^2(sM_{s-1}+s(s-1)M_{s-2}),
\end{align*}
so we have the first equation of statement.

To get the equations for $\M_s$ it is sufficient to use Lemmas \ref{MkMkMk}, \ref{M'kMkMk} and previous relation.
\epf


\section{Several Inequations}

Let denote by $\gg$ that the inequation $\ge$ holds for all coefficient before $x^n$ in the considering series. For example,
the notation $\sum a_nx^x\gg\sum b_nx^n$ means that for all $n$ the inequation $a_n\ge b_n$ holds. It is easy to verify that:

if $X\gg Y$ and $Z\gg 0$, then $XZ\gg YZ$;

if $X\gg Y$ and $V\gg W$, then $X+V\gg Y+W$.

\begin{lem}\label{otsenki} For $n>0$
$$
q^{n-1}\E(\nu_{n-1})^{\underline s} \le \E(\nu_n-1)^{\underline s}\le \E(\nu_{n-1})^{\underline s}
$$
\end{lem}
\pf Left inequation follows from:
$$
\M'_s=B'M_s\gg M_s
$$
with help of \hyperref[convolution]{convolution formula} and because of $B'\gg 1$. Right inequation follows from:
$$
\M'_s=B'M_s=B'AB^s=A'B^s=A(x/q)B^s\ll A(x/q)B(x/q)^s=M_s(x/q).
$$
\epf

\begin{lem}\label{otsenki2} For all $n\ge 1$ and $s\ge 1$ the following inequations hold:
$$
(n-1)^{\underline s}\cdot q^{(n-1)s} \le \E(\nu_n-1)^{\underline s}\le 2(n-1)^{\underline s}q^{(n-1)(s+1)/2}
$$
\end{lem}
\pf
Left inequation.
$$
\frac{x\M_s'}{s!}=\frac{xB'M_s}{s!}=xA'\frac{B^s}{s!}=xA'B_s,
$$
hence, by the \hyperref[convolution]{convolution formula} we get:
$$
\frac{n\E(\nu_n-1)^{\underline s}}{s!}=\sum_{k=0}^n\binom nk q^{k(n-k)}k\P\{\nu_{n-k}=s\},
$$
where the last summation we can estimate by the summand as $k=n-s$, and therefore we have:
$$
\E(\nu_n-1)^{\underline s}\ge \binom{n}{n-s}s! q^{s(n-s)}q^{s(s-1)/2}\cdot\frac{n-s}{n}=(n-1)^{\underline s}\cdot q^{(n-1)s}q^{-s(s-1)/2},
$$
here we get the left equation of Lemma statement.

Right inequation. Following relations one can get from the results that were proved above.
\begin{align*}
\frac{x(\widehat\M_s)'}{s!}&=x\left(\frac{x^2}{2q}(B')^2\frac{M_{s-1}}{(s-1)!}\right)'=
\frac{x^2(B')^2+x^3B'B''}{q}\frac{M_{s-1}}{(s-1)!}+\frac{x^3(B')^2}{2q}\frac{M_{s-1}'}{(s-1)!}=\\
&=\frac{x^2}{q}\left((B')^2\frac{M_{s-1}}{(s-1)!}+xB'B''\frac{M_{s-1}}{(s-1)!}+\frac x2(B')^3\frac{M_{s-1}}{(s-1)!}+\frac x2(B')^3\frac{M_{s-2}}{(s-2)!}\right)\\
\frac{x^2}{qs!}\M_s''&=\frac{x^2}{qs!}(B'M_s)'=\frac{x^2}{q}\left(B''\frac{M_s}{s!}+(B')^2\frac{M_s}{s!}+(B')^2\frac{M_{s-1}}{(s-1)!}\right)
\end{align*}

\begin{align*}
\frac{x(\widehat\M_s)'}{s!}-s\frac{x^2}{qs!}\M_s''&=
\frac{x^2}{q}(B')^2\frac{M_{s-1}}{(s-1)!}(1-s)
+\frac{x^2}{q}B''M_{s-1}\left(\frac{xB'}{(s-1)!}-\frac{sB}{s!}\right)\\
&+\frac{x^2}{q}(B')^2M_{s-1}\left(\frac{xB'}{2(s-1)!}-\frac{sB}{s!}\right)
+\frac{x^3}{2q}(B')^3\frac{M_{s-2}}{(s-2)!}
\end{align*}

\begin{align}
x(\widehat\M_s)'-\frac{sx^2}{q}\M_s''=&
s\frac{x^2}{q}B''M_{s-1}(xB'-B)+\frac{sx^2}{q}(B')^2M_{s-1}(xB'/2-B)\nonumber\\
&+s(s-1)\frac{x^2}{q}(B')^2M_{s-2}(xB'/2-B)=\nonumber\\
=&s\frac{x^2}{q}B''M_{s-1}(xB'-B)+\frac{s!x^2}{q}(B')^2\left(\frac{M_{s-1}}{(s-1)!}+\frac{M_{s-2}}{(s-2)!}\right)(xB'/2-B)=\nonumber\\
=&s\frac{x^2}{q}B''M_{s-1}(xB'-B)+2\widehat M_s(xB'/2-B).\label{Leq1}
\end{align}

Since $n-1\ge 0$, $n/2-1\ge 0$ for $n\ge 2$, $n/2-1\ge -1/2$ for $n=1$ it follows that
$$
xB'-B\gg 0;\qquad \frac{xB'}{2}-B\gg -\frac{x}{2},
$$
and from the equations \eqref{Leq1} we get the next inequation:
$$
x(\widehat\M_s)'+x\widehat M_s\gg \frac{sx^2}{q}\M_s''.
$$
Now we get coefficients before $x^n$:
$$
n(\E(\nu_n-1)^{\underline s})'_q+nq^{n-1}(\E(\nu_{n-1})^{\underline s})'_q\ge \frac{sn(n-1)}{q}\E(\nu_n-1)^{\underline s}.
$$
Dividing by n we get:
\begin{equation}\label{Leq2}
(\E(\nu_n-1)^{\underline s})'_q+q^{n-1}(\E(\nu_{n-1})^{\underline s})'_q\ge \frac{s(n-1)}{q}\E(\nu_n-1)^{\underline s}\quad\mbox{as } n>0.
\end{equation}
It is esy to see that
$$
q^{n-1}(\E(\nu_{n-1})^{\underline s})'_q=(q^{n-1}\E(\nu_{n-1})^{\underline s})'_q-\frac{(n-1)}{q}q^{n-1}\E(\nu_{n-1})^{\underline s},
$$
where we use derivative of product. Therefore from this and \eqref{Leq2} we get
\begin{equation}\label{predv}
(\E(\nu_n-1)^{\underline s}+q^{n-1}\E(\nu_{n-1})^{\underline s})'_q\ge \frac{s(n-1)}{q}\E(\nu_n-1)^{\underline s}+\frac{(n-1)}{q}q^{n-1}\E(\nu_{n-1})^{\underline s}
\end{equation}
Hence, dividing by $\E(\nu_n-1)^{\underline s}+q^{n-1}\E(\nu_{n-1})^{\underline s}$ we find the inequation
\begin{equation}\label{Leq3}
\frac{(\E(\nu_n-1)^{\underline s}+q^{n-1}\E(\nu_{n-1})^{\underline s})'_q}{\E(\nu_n-1)^{\underline s}+q^{n-1}\E(\nu_{n-1})^{\underline s}}
\ge\frac{n-1}{q}\cdot\frac{s\E(\nu_n-1)^{\underline s}+q^{n-1}\E(\nu_{n-1})^{\underline s}}{\E(\nu_n-1)^{\underline s}+q^{n-1}\E(\nu_{n-1})^{\underline s}}
\end{equation}
Note, that the function $f(t)=(s+t)/(1+t)$ not increases as $t$ increases, if $s\ge 1$ and $t>0$. From Lemma \ref{otsenki} it follows that:
$$
\E(\nu_n-1)^{\underline s}\ge q^{n-1}\E(\nu_{n-1})^{\underline s}\quad\mbox{as } n>0
$$
Therefore from this and \eqref{Leq3} we get:
\begin{equation}\label{n_n-1_usilenie}
\frac{(\E(\nu_n-1)^{\underline s}+q^{n-1}\E(\nu_{n-1})^{\underline s})'_q}{\E(\nu_n-1)^{\underline s}+q^{n-1}\E(\nu_{n-1})^{\underline s}}\ge\frac{(s+1)(n-1)}{2q}
\end{equation}

Let $q\le q_1\le 1$, and let $\E_1=\E|_{q=q_1}$. By integrating \eqref{n_n-1_usilenie} on the interval $[q;q_1]$ we get follows:
$$
\ln\left(\E(\nu_n-1)^{\underline s}+q^{n-1}\E(\nu_{n-1})^{\underline s}\right)\Big|_q^{q_1}\ge 
\frac{(s+1)(n-1)}{2}\ln q\Big|_q^{q_1}
$$
Then we put both sides of this inequation into the argument of function $e^x$, and get:
$$
\E(\nu_n-1)^{\underline s}+q^{n-1}\E(\nu_{n-1})^{\underline s}\le 
\left(\E_1(\nu_n-1)^{\underline s}+q_1^{n-1}\E_1(\nu_{n-1})^{\underline s}\right)\left(\frac{q^{n-1}}{q_1^{n-1}}\right)^{(s+1)/2}
$$
or
$$
\E(\nu_n-1)^{\underline s}\le 
\left(\E_1(\nu_n-1)^{\underline s}+q_1^{n-1}\E_1(\nu_{n-1})^{\underline s}\right)\left(\frac{q^{n-1}}{q_1^{n-1}}\right)^{(s+1)/2}.
$$
Then we set $q_1=1$ and finally find that
$$
\E(\nu_n-1)^{\underline s}\le 2(n-1)^{\underline s}q^{(n-1)(s+1)/2},
$$
because as $q_1=1$ we have $\E_1(\nu_n-1)^{\underline s}=(n-1)^{\underline s}$.

This proofs the right inequation of Lemma.
\epf

Note, that if we put $s=n-1$, then we have an equation
$$
\E(\nu_n-1)^{\underline{n-1}}=(n-1)!q^{n(n-1)/2}=(n-1)^{\underline{n-1}}q^{n(n-1)/2},
$$
where the right hand side is equal to half of the just proved estimation.


\section{Asymptotic behavior of $\nu_n$}

In this section we consider an asymptotics of moments $\E(\nu_n-1)^{\underline s}$ as $s$ is fixed and positive.
We will study a behavior of moments in the following zones of parameters:
\begin{enumerate}
\item $p\to 0$, $n=\const$;
\item $q^n\to e^{-\al}$, where fixed $\al\ge 0$ and $n\to\infty$;
\item $q^n\to 0$ as $n\to\infty$
\begin{enumerate}
\item[3.1] $nq^n\to\infty$,
\item[3.2] $nq^n\to\al>0$,
\item[3.3] $nq^n\to 0$ (in this case $p$ can be a positive constant $<1$).
\end{enumerate}
\end{enumerate}

\subsection{Asymptotics for $n=\const$}

It is easy to prove the following Lemma, besause the minimal graph with $n$ vertices and $s$ components is a forest with $s$ trees.
\begin{lem}\label{n-const}
If $p\to 0$ and $n=\const$, then for any $s\le n$ the following equation holds:
$\P\{\nu_n=s\}=F_{s,n}p^{n-s}+O(p^{n-s+1})$. In particular, $p_n=n^{n-2}p^{n-1}+O(p^n)$.
\end{lem}
Hence we have the following
\begin{thrm}\label{zone1}
If $p\to 0$ and $n=\const$, then for any $s\le n$:
$$
\E\nu_n^s=n^s(1+o(1)),\quad \E(\nu_n-1)^{\underline s}\to n^{\underline s}.
$$
\end{thrm}

\subsection{Asymptotics for $q^n\to e^{-\al}$}

Let
$$
\be(x)=\sum_{k=1}^\infty\frac{k^{k-2}}{k!}x^k.
$$
This series converges as $|x|\le e^{-1}$ and this is a generated function for sequence of numbers of labelled trees \cite{Kolchin}.

\begin{thrm}\label{np-const}
If $q^n\to e^{-\al}$ as $n\to\infty$, where $\al\ge 0$, then for $s\ge 0$
$$
\E(\nu_n-1)^{\underline s}= \left(\frac n\al\be(\al e^{-\al})\right)^s(1+o(1)).
$$
In particular, for $\al=0$ we have the relation $\E(\nu_n-1)^{\underline s}\sim n^s$.
\end{thrm}
\pf We will use a mathematical induction on the parameter $s$. It is clear that the statement of Theorem holds for $s=0$.
Let us suppose that it holds for $s-1$ and will show it for $s\ge 1$.

From the relation $x\M'_s=xB'M_s=xB'BM_{s-1}=Bx\M'_{s-1}$ (see Lemma \ref{M'kMkMk}) and from the \hyperref[convolution]{convolution formula} we get:
\begin{multline*}
n\E(\nu_n-1)^{\underline s}=\sum_{k=0}^n\binom nkq^{k(n-k)}p_k(n-k)\E(\nu_{n-k}-1)^{\underline {s-1}}=\\
=n\sum_{k=0}^{n-1}\binom {n-1}kq^{k(n-k)}p_k\E(\nu_{n-k}-1)^{\underline {s-1}}=n(S_1+S_2),
\end{multline*}
where
\begin{align*}
S_1=&\sum_{k=0}^{k_0}\binom {n-1}kq^{k(n-k)}p_k\E(\nu_{n-k}-1)^{\underline {s-1}},\\
S_2=&\sum_{k=k_0+1}^{n-1}\binom {n-1}kq^{k(n-k)}p_k\E(\nu_{n-k}-1)^{\underline {s-1}}
\end{align*}
As $k$ is fixed, one can get next relations: $\binom {n-1}kq^{k(n-k)}\sim (nq^n)^k/k!$, $p_k\sim k^{k-2}p^{k-1}$ (Lemma \ref{n-const}).
And from the induction hypothesis we get: $\E(\nu_{n-k})^{\underline {s-1}}\sim\left(\frac n\al\be(\al e^{-\al}\right)^{s-1}$. Hence:

$$
S_1=\sum_{k=0}^{k_0}\frac{(npq^n)^k}{pk!}k^{k-2}\left(\frac n\al\be(\al e^{-\al})\right)^{s-1}(1+o(1))=
\sum_{k=0}^{k_0}\frac{n}{\al}\frac{k^{k-2}}{k!}(\al e^{-\al})^k\left(\frac n\al\be(\al e^{-\al})\right)^{s-1}(1+o(1)),
$$
where we use the asymptotics $np\to \al$ and $npq^n\to\al e^{-\al}$, which is follows from Theorem conditions.

So, it is easy to see that $S_1/\left(\frac n\al\be(\al e^{-\al})\right)^s$ as closed to 1 as $k_0$ is bigger, because of the convergence of the series $\be(x)$ for $x=\al e^{-\al}$.

Let we estimate $S_2$. It is clear that $\E(\nu_{n-k}-1)^{\underline {s-1}}\le n^{s-1}$. From this and from the equation \eqref{nkpk} we get
\begin{multline*}
1\ge \sum_{k=k_0+1}^{n-1}\binom {n-1}kq^{k(n-1-k)}p_k\frac{k}{n-1}\ge\frac{k_0}{n-1}\sum_{k=k_0+1}^n\binom {n-1}kq^{k(n-k)}p_k\ge\\
\ge\frac{k_0}{n^s}\sum_{k=k_0+1}^{n-1}\binom {n-1}kq^{k(n-k)}p_k\E(\nu_{n-k}-1)^{\underline {s-1}}=\frac{k_0}{n^s}S_2.
\end{multline*}
Therefore,
$$
S_2=\frac{1}{k_0}O\left(n^s\right)=\frac{1}{k_0}O\left(\frac n\al\be(\al e^{-\al})\right)^s,
$$
i.e. the ratio $S_2/\left(\frac n\al\be(\al e^{-\al})\right)^s$ tends to 0 as $k_0\to\infty$.

Thus, $\E(\nu_n-1)^{\underline s}= \left(\frac n\al\be(\al e^{-\al})\right)^s(1+o(1))$.

In the case of $\al=0$ the proof of Theorem is similary, but instead of $\be(\al e^{-\al})/\al$ we should write 1 at all places.
\epf

\subsection{Asymptotics for $q^n\to 0$}

\begin{thrm}\label{middle} Let $q^n\to 0$ and $nq^n\ge C$ as $n\to\infty$, where fixed $C>0$, then
$$
\E(\nu_n-1)^{\underline s}= (nq^n)^s(1+o(1)).
$$
\end{thrm}
\pf We will use an induction by $s$. It is clear that the statement of Theorem holds for $s=0$.
Let us suppose that it holds for $s-1$ and will show it for $s\ge 1$.

The following relations hold:
\begin{multline*}
\sum_{n=0}^\infty \frac{(\E(\nu_n)^{\underline s}-\E(\nu_{n+1}-1)^{\underline s})x^n}{n!q^{n(n-1)/2}}=M_s-\M'_s(xq)=
M_s-B'(xq)M_s(xq)=\\
=AB^s-A'(xq)B^s(xq)=A(B^s-B^s(xq))=A(s!B_s-s!B_s(xq))=\\
=A\sum_{n=0}^\infty (1-q^n)\frac{x^n\P\{\nu_n=s\}s!}{n!q^{n(n-1)/2}}
\ll A\sum_{n=0}^\infty np\frac{x^n\P\{\nu_n=s\}s!}{n!q^{n(n-1)/2}}=\\
=Apx(B^s)'=spxAB'B^{s-1}=spxB'M_{s-1}=spx\M'_{s-1}=sp\sum_{n=0}^\infty\frac{\E(\nu_n-1)^{\underline {s-1}}nx^n}{n!q^{n(n-1)/2}},
\end{multline*}
where we use the fact, that $(1-q^n)\le n(1-q)=np$. Therefore we get that
$$
\E(\nu_n)^{\underline s}-\E(\nu_{n+1}-1)^{\underline s}\le spn\E(\nu_n-1)^{\underline {s-1}}
$$
or
$$
\E(\nu_{n-1})^{\underline s}\le \E(\nu_n-1)^{\underline s}+spn\E(\nu_{n-1}-1)^{\underline {s-1}}.
$$
From the equation \eqref{Leq3} it follows that
\begin{multline}\label{drob}
\frac{(\E(\nu_n-1)^{\underline s}+q^{n-1}\E(\nu_{n-1})^{\underline s})'_q}{\E(\nu_n-1)^{\underline s}+q^{n-1}\E(\nu_{n-1})^{\underline s}}\ge
\frac{n-1}{q}\cdot\frac{s\E(\nu_n-1)^{\underline s}+q^{n-1}\E(\nu_{n-1})^{\underline s}}{\E(\nu_n-1)^{\underline s}+q^{n-1}\E(\nu_{n-1})^{\underline s}}\ge \\
\frac{n-1}{q}\cdot\frac{s\E(\nu_n-1)^{\underline s}+q^{n-1}(\E(\nu_n-1)^{\underline s}+snp\E(\nu_{n-1}-1)^{\underline {s-1}})}{\E(\nu_n-1)^{\underline s}+q^{n-1}(\E(\nu_n-1)^{\underline s}+snp\E(\nu_{n-1}-1)^{\underline {s-1}})}=\\
=\frac{n-1}{q}\cdot\frac{s+q^{n-1}+snpq^{n-1}\E(\nu_{n-1}-1)^{\underline {s-1}}/\E(\nu_n-1)^{\underline s}}{1+q^{n-1}+snpq^{n-1}\E(\nu_{n-1}-1)^{\underline {s-1}}/\E(\nu_n-1)^{\underline s}}.
\end{multline}
By the induction hypothesis and from Lemma \ref{otsenki2} we get that
$$
\frac{\E(\nu_{n-1}-1)^{\underline {s-1}}}{\E(\nu_n-1)^{\underline s}}\le C_1\frac{(nq^n)^{s-1}}{(n-1)^{\underline s}q^{(n-1)s}}\le \frac{C_2}{nq^n},
$$
where the positive constants $C_k$, generally speaking, are depends on the parameter $s$. By putting this inequation into \eqref{drob} we get:
\begin{multline}\label{ner-vo}
\frac{(\E(\nu_n-1)^{\underline s}+q^{n-1}\E(\nu_{n-1})^{\underline s})'_q}{\E(\nu_n-1)^{\underline s}+q^{n-1}\E(\nu_{n-1})^{\underline s}}\ge
\frac{n-1}{q}\cdot\frac{s+q^{n-1}+C_4npq^{n-1}/(nq^n)}{1+q^{n-1}+C_4npq^{n-1}/(nq^n)}\ge\\
\ge\frac{n-1}{q}(s-sq^{n-1}-C_5npq^{n-1}/(nq^n))\ge\frac{n-1}{q}(s-sq^{n-1}-C_6p)\ge\\
\ge (n-1)sq^{-1}-C_7(n-1)q^{n-2}-C_8np.
\end{multline}

Let $q_1=\ep^{1/(n-1)}$, where $\ep$ is an arbitrary small positive number, hence $q_1^{n-1}=\ep$ and $q<q_1$ (it follows from $q^n\to 0$).
Besides let denote $\E_1=\E|_{q=q_1}$ as it was above.

Now, we integrate the inequation \eqref{ner-vo} on the interval $[q;q_1]$ and get that

$$
\ln\left(\E(\nu_n-1)^{\underline s}+q^{n-1}\E(\nu_{n-1})^{\underline s}\right)\Big|_q^{q_1}\ge 
s(n-1)\ln q\Big|_q^{q_1}-(n-1)C_7\frac{q^{n-1}}{n-1}\Big|_q^{q_1}-C_8n(q-q^2/2)\Big|_q^{q_1}
$$
or
\begin{multline}\label{q-q1}
\E(\nu_n-1)^{\underline s}+q^{n-1}\E(\nu_{n-1})^{\underline s}\le\\
\le \left(\E_1(\nu_n-1)^{\underline s}+q_1^{n-1}\E_1(\nu_{n-1})^{\underline s}\right)
\left(\frac{q}{q_1}\right)^{s(n-1)}e^{C_7(q_1^{n-1}-q^{n-1})}e^{C_8n(q_1-q+q^2/2-q_1^2/2)}\le\\
\le\frac{\E_1(\nu_n-1)^{\underline s}+q_1^{n-1}\E_1(\nu_{n-1})^{\underline s}}{\ep^s}q^{s(n-1)}e^{C_9\ep},
\end{multline}
because $q_1^{n-1}=\ep$ and $n(q_1-q+q^2/2-q_1^2/2)=n(q_1-q)(1-q/2-q_1/2)=n(p-p_1)(p/2+p_1/2)\le np^2 \to 0$.
The last expression is follows from $np^2\cdot nq^n=(np)^2e^{n\ln q}\to 0$ and from the conditions of Theorem.

By Theorem \ref{np-const} we get that
$$
\E_1(\nu_n-1)^{\underline s}=\left(\frac n\al\be(\al e^{-\al})\right)^s(1+o(1))=n^s\be(\al \ep)^s/\al^s(1+o(1)).
$$
where $\al=-\ln\ep$.

Besides that,
\begin{multline*}
q_1^{n-1}\E_1(\nu_{n-1})^{\underline{s}}=\ep(\E_1(\nu_{n-1}-1)^{\underline{s}}+s\E_1(\nu_{n-1}-1)^{\underline{s-1}})=\\
=\ep(n^s\be(\al \ep)^s/\al^s+sn^{s-1}\be(\al \ep)^{s-1}/\al^{s-1})(1+o(1)),
\end{multline*}
because $q_1^{n-1}=\ep$ and again from Theorem \ref{np-const}. From this and from \eqref{q-q1} it follows that
$$
\E(\nu_n-1)^{\underline s}\le\frac{(1+\ep)\be(\al \ep)^s/\al^s+\ep s\be(\al \ep)^{s-1}/\al^{s-1}}{\ep^s} n^sq^{sn}e^{C_{10}\ep}(1+o(1)).
$$
Therefore, by choosing an arbitrary small $\ep>0$ and using the relationship $\be(x)\sim x$ as $x\to 0$ we get the relation:
$$
\limsup_{n\to\infty}\frac{\E(\nu_n-1)^{\underline s}}{(nq^n)^s}\le \lim_{\ep\to 0}\frac{(1+\ep)\be(\al \ep)^s+s\ep\al\be(\al \ep)^{s-1}}{\al^s\ep^s}e^{C_{10}\ep}=1.
$$
From Lemma \ref{otsenki2} we have $\E(\nu_n-1)^{\underline s}\ge (n-1)^{\underline s}q^{s(n-1)}=(nq^n)^s(1+o(1))$. Now we see that Theorem follows from theese both equations.
\epf

Let $nq^n\to\al$, where $\al$ is a positive constant. From Theorem \ref{middle} we see that $\E(\nu_n-1)^{\underline s}\to\al^s$.

It is known that in this case random variable $(\nu_n-1)$ tends to Poisson distribution with the parameter $\al$.

Thus we have
\begin{thrm}\label{poisson}
If $nq^n\to\al$ as $n\to\infty$ and $\al$ is a fixed positive constant, then for any fixed integer $k\ge 1$:
$$
\P\{\nu_n=k\}\to \frac{\al^{k-1}}{(k-1)!}e^{-\al}.
$$
\end{thrm}

From Lemma \ref{otsenki2} it follows that if $nq^n\to 0$, then $\E(\nu_n-1)\asymp nq^{n-1}$. So we can conclude that $\nu_n$ tends to 1.
Below we'll show an estimation of $p_n$ in this case.

\section{Several Consequences}

Generally, we can conclude that in all zones of parameters $p$ and $n$
$$
\E(\nu_n-1)^{\underline s} \sim (\be(npq^n)/p)^s\quad\mbox{ as } nq^n\to\infty
$$
and
$$
\E(\nu_n-1)^{\underline s} = (\be(npq^n)/p)^s +o(1)\quad\mbox{ as } nq^n=O(1)
$$

It is easy to verify, because if $np\to\infty$ or $np\to 0$, then it follows that $npq^n\to 0$ and $\be(npq^n)/p\sim nq^n$.

Now we can estimate the probability $p_n$ that graph $G(n,p)$ is connected.
\begin{equation}\label{Leq7}
p_n=\P\{\nu_n<2-1/n\}=1-\P\{\nu_n-1\ge 1-1/n\}\ge 1-\E(\nu_n-1)\frac{n}{n-1},
\end{equation}
and from Lemma \ref{otsenki2} we get:
\begin{equation}\label{pn-big}
p_n\ge 1-2nq^{n-1}.
\end{equation}

If we put $p=\frac{c\ln n}{n}$ and $c>1$, then we have $nq^{n-1}=n\exp\{-c\ln n + O(\ln^2n)/n\}=n^{1-c}(1+O(\ln^2 n)/n)$.
Therefore we finally get:
\begin{equation}
p_n\ge 1-\frac{2}{n^{c-1}}(1+O(\ln^2 n)/n).
\end{equation}

If $nq^n\to\al$ (for example, $p=(\ln n + c +o(1))/n$, where $\al=e^{-c}$), then from Theorem \ref{poisson} we get that:
$$
p_n\to e^{-\al}.
$$

To estimate $p_n$ as $nq^n\to\infty$ we now consider the isolating probability.
Let $pi_n$ be a probability that $G(n,p)$ has an isolated vertex. Let $A_i$ be an event that $i$-th vertex is isolated, then from the Inclusion--exclusion principle we get:
$$
pi_n=\P\{A_1\cup\dots\cup A_n\}=\sum\limits_{k=1}^n(-1)^{k-1}\sum\limits_{1\le i_1<\dots<i_k\le n}\P\{A_{i_1}\dots A_{i_k}\}=
\sum\limits_{k=1}^n(-1)^{k-1}\binom nk\P\{A_1\dots A_k\}.
$$
It is easy to see that $\P\{A_1\dots A_k\}=q^{k(k-1)/2}q^{k(n-k)}$, so
$$
pi_n=\sum\limits_{k=0}^n(-1)^{k-1}\binom nkq^{k(n-k)}q^{k(k-1)/2}+1.
$$
According to \hyperref[convolution]{convolution formula} we can find that SG-function $PI$ of $\{pi_n\}$ is equal to $RT+A$, where $R$ and $T$ are SG-functions of the corresponding sequences $\{r_n\}$ and $\{t_n\}$, which are defined as follows:  $r_n=(-1)^{n-1}q^{n(n-1)/2}$ and $t_n=1$.

Hence we have
$$
R=\sum\limits_{n=0}^\infty \frac{r_nx^n}{n!q^{n(n-1)/2}}=-e^{-x};\qquad T=A.
$$
Thus $PI=A-e^{-x}A=A(1-e^{-x})$.

Since $(1-e^{-x})\le x$ it follows that $PI\ll Ax$, and from the \hyperref[convolution]{convolution formula} we obtain
\begin{equation}\label{pi-little}
pi_n\le nq^{n-1}.
\end{equation}

It is easy to see that $PI'=A'(1-e^{-x})+Ae^{-x}= PI\cdot B' + A-PI \gg PI\cdot B'$, because $A-PI\gg 0$, and from the \hyperref[convolution]{convolution formula} we get:
$$
n pi_n\ge\sum_{k=1}^n\binom nkq^{k(n-k)}pi_{n-k} k p_k\ge n(n-1)q^{n-1}p_{n-1}
$$
or
\begin{equation}\label{pi-pn}
pi_{n+1}\ge nq^np_n
\end{equation}

So, if $nq^n\to\al>0$, then $pi_n\ge \al e^{-\al} + o(1)$.

And also we have
\begin{equation}\label{pn-little}
p_n\le pi_{n+1}/(nq^n)\le 1/(nq^n)
\end{equation}

Since $PI=A(1-e^{-x})$ it follows that $PI e^x = Ae^x-A$ and, therefore, $(PI-A)(e^x-1)=-PI$. From the relation $e^x-1>x$ we get that
$PI\gg (A-PI)x$, therefore from the \hyperref[convolution]{convolution formula} we find that
$pi_n \ge (1-pi_n)nq^{n-1}$, then $(1-pi_n)\le 1/(nq^{n-1})$ and we get finally
\begin{equation}\label{pi-big}
pi_n\ge 1-\frac{1}{nq^{n-1}}
\end{equation}

Now we can combine all obtained results \eqref{pn-big}, \eqref{pn-little}, \eqref{pi-pn}, \eqref{pi-little} and \eqref{pi-big} in the following
\begin{thrm} For all $n\ge 1$
\begin{align*}
1-2nq^{n-1}&\le p_n\le\frac{1}{nq^n},\\
1-\frac{1}{nq^n}&\le pi_n\le nq^{n-1},\\
nq^np_n&\le pi_{n+1}
\end{align*}
And if $nq^n\ge C>0$ as $n\to\infty$, then we can substitute $nq^n$ by $\E(\nu_n-1)(1+o(1))$ in theese relations.
\end{thrm}


\begin{thebibliography}{99}
\bibitem{Erdos-Renyi} Erd\"os, P. and R\'enyi, A. (1959). "On Random Graphs." {\itshape Publicationes Mathematicae} {\bfseries 6}: 290-297.
\bibitem{Erdos-Renyi2} Erd\"os, P. and R\'enyi, A. (1960) "On the Evolution of Random Graphs." {\itshape Publ. Math. Inst. Hungar. Acad. Sci.} {\bfseries 5}, 17-61.
\bibitem{Bollobas} Bollob\'as, B. (2001) {\itshape Random Graphs (2nd ed.).} Cambridge University Press.
\bibitem{Kolchin} Kolchin, V. F. {\itshape Random Graphs.} New York: Cambridge University Press, 1998.
\bibitem{A062734} Sloane, N. J. A. Sequence \href{http://oeis.org/A062734}{A062734} in "The On-Line Encyclopedia of Integer Sequences."
\bibitem{Harari-Palmer} Harary, Frank; Palmer, Edgar M. (1973). {\itshape Graphical Enumeration.}
\end{thebibliography}
\end{document}